\documentclass[a4paper, 12pt]{article}
\usepackage{enumerate, theorem}
\usepackage{amsmath, amsfonts, amssymb}
\usepackage[height=22.5cm, width=15cm]{geometry}
\usepackage[all]{xy}
\usepackage{mathrsfs, yfonts}

\DeclareMathOperator{\C}{\mathbb{C}}

\newcommand{\parag}[1]{\paragraph{\sc{#1.}} }

\newtheorem{thm}{Theorem}[subsection]
\newtheorem{defn}[thm]{Definition}
\newtheorem{cor}[thm]{Corollary}
\newtheorem{prop}[thm]{Proposition}
\newtheorem{lemma}[thm]{Lemma}

\begin{document}

\title{Gauduchon's form and compactness of the space of divisors (second version)}

\date{16/05/17}

 \author{Daniel Barlet\footnote{Institut Elie Cartan, Alg\`{e}bre et  G\'eom\`{e}trie,\newline
Universit\'e de Lorraine, CNRS UMR 7502   and  Institut Universitaire de France.}.}

 \maketitle
 
 $ \hfill ${\it En hommage \`{a} E. Bishop }
 
 \parag{Abstract} We show that in a holomorphic family of compact complex connected manifolds parametrized by an irreducible complex space $S$, assuming that on  a dense Zariski open set $S^{*}$ in $S$ the fibres satisfy the $\partial\bar\partial-$lemma, the algebraic dimension of each fibre in this family is at least equal to the minimal algebraic dimension of the fibres in $S^{*}$. For instance, if each fibre in $S^{*}$ are Moishezon, then all fibres are Moishezon.\\
 
 \parag{AMS Classification 2010} 32 G 10 - 32 J 18- 32 J 27- 32 Q 99- 32 G 05- 32 J 10.\\
 
 \parag{Key words}Family of compact complex manifods,  Algebraic dimension, relative codimension 1 cycle-space, $\partial\bar\partial-$lemma, strongly Gauduchon manifolds.\\
 
  \section{Introduction}
 
 In this  article we give a rather elementary proof of the following result which answers a ``classical'' question.
 
 \begin{thm}\label{main}
 Let $\pi : \mathcal{X} \to D$ be a smooth holomorphic family of compact complex connected manifolds parametrized by the unit disc $D$ in $\C$. Assume  that for each $s \in D^{*}$ the fiber $X_{s}$ of $\pi$ at $s$ satisfies the $\partial\bar\partial-$lemma. Let
  $$a := \inf\{ alg(X_{s}), s \in D^{*} \}$$
  where $alg(X)$ denotes the algebraic dimension of the compact complex connected manifold $X$. Then  we have $alg(X_{0}) \geq a $.
 \end{thm}
 
 As a special case (when $a = n := dim X_{s} $) we obtain that, if for any $s \in D^{*}$ each $X_{s}$ is a Moishezon manifold,  $X_{0}$ is also a Moishezon manifold.
 
 The previous theorem gives easily the following corollary (see [B.15] for details).
 
 \begin{cor}\label{ variant}
 Let $\pi : \mathcal{X} \to S$ be a holomorphic family of compact complex connected manifolds parametrized by a reduced and irreducible complex space $S$. Let $S^{*}$ be a dense Zariski open set in $S$.  Assume that  for each $s \in S^{*}$ the fiber $X_{s}$ of $\pi$ at $s$ satisfies the $\partial\bar\partial-$lemma. Let $a := \inf\{ alg(X_{s}), s \in S^{*} \}$. Then for any $s \in S$ we have $alg(X_{s}) \geq a $.$\hfill \blacksquare$
 \end{cor}
 
 Note that  the proof  in [B.15] shows that the minimum of the algebraic dimension is obtained at the general point\footnote{this means on the complement of a countable union of closed nowhere dense analytic subsets in $S$ (see the proposition \ref{classique} below).} in $S$.\\
 
 In the first section we show that the existence of a Gauduchon metric on a compact complex connected manifold $X$ implies the compactness of the connected components of the space of divisors in $X$. This new proof of this classical result is the key of our proof for the theorem above which appears as a relative version of it.\\
 
 {\it I thank J-P. Demailly for some constructive comments which help to improve this article.}

 \section{The absolute case}
 
 Let me begin by two simple lemmas.
 
 \begin{lemma}\label{minimum}
 Let $M$ be a connected reduced complex space and let $T : M \to \mathbb{R}^{+}$ be a continuous function on $M$. Assume that $T$ is pluri-harmonic on the smooth part $M'$ of $M$ and that the function $T$ achieves its minimum at a point $x_{0}\in M$. Then $T$ is constant on $M$.
 \end{lemma}
 
 \parag{Proof} This is an elementary exercice.$\hfill \blacksquare$
 
 \begin{lemma}\label{volume function}
 Let $X$ be a compact  reduced complex space and let $\omega$ be a continuous  real $(q, q)-$form on $X$ which is strictly positive in the Lelong sense (see section 3.1). Let $\Gamma$ be a connected component of the reduced complex  space of compact $q-$cycles in $X$ and define
 $$ \theta : \Gamma \to \mathbb{R}^{+} \quad {\rm by}  \quad \gamma \mapsto \theta(\gamma) : = \int_{\gamma} \omega .$$
 Then the continuous function $\theta$ achieves its minimum on $\Gamma$.
 \end{lemma}
 
 \parag{Proof} Note first that $\theta$ is continuous thanks to the prop. IV 3.2.1 of [B-M 1]. Let $\alpha := \inf \{ \theta(\gamma), \gamma \in \Gamma \}$. Then the subset $A := \{\gamma \in \Gamma \ / \ \theta(\gamma) \leq \alpha +1\}$ is a compact subset in $\Gamma$ thanks to Bishop's theorem [Bi.64], because it is a closed subset in $\Gamma$ such that each cycle in $A$ has bounded volume (relative to $\omega$) and support in the compact space $X$ (see th. IV 2.7.23 in [B-M 1]). Then $\theta$ achieves its minimum on $A$.$\hfill \blacksquare$
 
 \begin{prop}\label{application}
  Let $X$ be a compact  reduced complex space and let $\omega$ be a continuous   real $(q, q)-$form on $X$ which is strictly positive in the Lelong sense (see section 3.1) and satisfies $\partial\bar\partial \omega = 0$ as a current  on $X$. Then any connected component of the reduced complex space $\mathcal{C}_{q}(X)$ of compact $q-$cycles in $X$ is compact.
  \end{prop}
  
  \parag{Proof} The function $\gamma \to \theta(\gamma) :=  \int_{\gamma} \omega$ is continuous on $\mathcal{C}_{q}(X)$ and $\partial\bar\partial-$closed in the sense of currents;  as its achieves its minimum on any connected component  $\Gamma$  of $\mathcal{C}_{q}(X)$ thanks to the lemma \ref{volume function} , we conclude that it is constant on each  $\Gamma$,  thanks to lemma \ref{minimum}, and then any such  $\Gamma$ is compact (see [B-M 1] ch.IV   th. 2.7.20).$\hfill\blacksquare$
  
  \parag{Remark} As any compact complex manifold $X$ of dimension $n$ admits a Gauduchon metric (see [G.77]), the proposition above gives a proof of the fact that the space of divisors of such a $X$ has always compact connected components, which is a well known classical result (see [C.82] and [Fu.82]).

  \section{The relative case}
  
 \subsection{Relative strong positivity in the Lelong sense}
 
 Let $ \pi: X  \to S $ be a surjective holomorphic map between two irreducible complex spaces. Let $T_{\mathcal{X}/D}$ be the tangent Zariski linear space of $\pi$, so the kernel of the tangent map 
 $$ T_{\pi} : T_{\mathcal{X}} \to \pi^{*}(T_{D}).$$
  Let $1 \leq q \leq \dim X - \dim S := n$ \ and let $\omega$ a continuous relative $(q, q)-$differential form on $X$. Then  $\omega$ induces  a continuous hermitian form of the fibres on the linear space   on $X$  associated to the coherent sheaf $\Omega_{X/S}^{q}:= \Omega_{X}^{q}\big/\pi^{*}(\Omega^{1}_{S})\wedge \Omega_{X}^{q-1}$.\\
   Let $p : Gr_{q}(X/S) \to X$ be the grassmannian of $q-$planes  in $T_{X/S}$, $U $ the universal $q-$vector bundle on $Gr_{q}(X/S)$ and $\theta : \Lambda^{q}(U) \to Gr_{q}(X/S)$ the line bundle which is the determinant of $U$. A $\pi-$relative continuous $(q, q)-$form on $X$ defines a continuous  hermitian form on $\Lambda^{q}(U)$.
 
 \begin{defn}\label{relative Lelong}
 We shall say that  a $\pi-$relative continuous $(q, q)-$form $\omega$  on  $X$ is {\bf strongly Lelong positive} at $x_{0}$ if the hermitian form on $\Lambda^{q}(U)$ defined by $\omega$ is a positive hermitian form at each point of $p^{-1}(x_{0})$ (so a continuous hermitian metric on this line bundle) .
 \end{defn}
 
 In other word that means that for any $q-$plane $P$ in the Zariski tangent space $T_{X, x_{0}}$ which is vertical (i.e. contained in the kernel of $T_{\pi, x_{0}}$) then $\omega_{x_{0}}[v_{1}\wedge \dots \wedge v_{q}] >  0$ when $v_{1}, \dots, v_{q}$ is a basis of $P$.\\
 As the map $p :  Gr_{q}(X/S) \to X$ is proper, the condition above is open. Moreover, it is an easy exercice to prove the following properties:
 \begin{enumerate}
 \item Assume that the  $\pi-$relative continuous $(q, q)-$form $\omega$  on  $X$  is strongly Lelong positive at $x_{0}$. Then  for any continuous hermitian metric $h$ on $X$ there exists a neighbourhood $V$ of $x_{0}$ in $X$ and a positive constant $C(V)$ such that on the open set $p^{-1}(V)$ the metric induced by $h^{\wedge q}$ and $\omega$ on the line bundle $\Lambda^{q}(U)$ satisfy
 $$ h^{\wedge q} \leq C(V). \omega .$$
 \item An easy consequence of the above estimate is the fact  that on any closed (complex) analytic subset $Z$ of pure dimension $q$ in a fiber of $\pi$ over $V$ we have the inequality of volume forms on $Z$
 $$  h^{\wedge q}_{\vert Z} \leq C(V).\omega_{\vert Z} .$$
 \item Assume now that $\omega$  is strongly Lelong positive on an open set $V$ in $X$ containing a compact  set $L$. Then there exists a positive constant $C_{L}$ such that for any compact relative $q-$cycle $\gamma$ contained in $L$ we have the inequality 
 $$ \int_{\gamma} \ h^{\wedge q} \leq C_{L}.\int_{\gamma} \ \omega.$$
 \end{enumerate}
 
 This means that the integral over compact relative $q-$cycles of a continuous $\pi-$relative $(q, q)-$form which is strongly Lelong positive on $X$  controls the volume of these cycles contained in a given  compact set  $L \subset X$ for any given hermitian metric on $X$.\\
 
 The following result follows easily of the descrition of compact subsets in the space of compact  cycles which is a consequence of E. Bishop's theorem.
 
 \begin{thm}\label{propre}
 Let $\pi : \mathcal{X} \to S$ a proper surjective map between irreducible complex spaces. Assume that there exists a smooth $d-$closed $2q-$differential form $\tilde{\omega}$ on $\mathcal{X}$ such that its $(q, q)-$part induces a  relative strongly Lelong positive form on $X$. Then any connected component of $\mathcal{C}_{q}(\pi)$, the space of relative  compact $n-$cycles of $\pi$, is proper over $S$.
 \end{thm}
 
 \parag{proof} Take such a connected component $\Gamma$. Then the continuous function given by  $\gamma \mapsto \int_{\gamma} \tilde{\omega}$ is constant on $\Gamma$. This implies that for any compact subset $K$ in $S$  the volume of cycles in $\Gamma \cap\pi^{-1}(K)$ for any continuous metric is bounded. But then, as $\pi^{-1}(K)$ is compact in $X$ this implies that $\Gamma \cap\pi^{-1}(K)$ is a compact subset in $\mathcal{C}_{q}(\pi) \subset \mathcal{C}_{q}(X)$ (see [B-M 1] IV th.2.7.20), concluding the proof.$\hfill\blacksquare$\\
 
\parag{Remark} In the case of smooth family of compact connected manifolds  on a complex disc (or polydisc) $\pi : \mathcal{X} \to D$ each of them satisfying the $\partial\bar\partial-$lemma, 
 the existence of a smooth $(q, q)-$relative differential form $\omega$ such that $\partial_{/\pi}\bar\partial_{/\pi} \omega = 0$ on $\mathcal{X}$ which is strongly Lelong positive on the fibres of $\pi$ allows to prove, in the same way than in the case of relative divisors, the properness over $D$ of the connected component of the relative cycle's space $\mathcal{C}_{q}(\pi)$ :\\
 the first step is to produce for any given $s \in D$, using the $\partial\bar\partial-$lemma on $X_{s}:= \pi^{-1}(s)$ smooth forms $\alpha$ and $\beta$ on $X_{s}$ of type $(q, q-1)$ and $(q-1, q)$ respectively such that the form $\omega -\partial\alpha - \bar\partial\beta$ is $d-$closed on $X_{s}$, but with $(q, q)-$ part strongly Lelong positive on $X_{s}$. Then, using the local $\mathscr{C}^{\infty}$ triviality of $\pi$ around $s$ to produce on $\pi^{-1}(\Delta_{s})$, where $\Delta_{s}$ is a small open (poly-)disc around $s$ in $D$, a smooth $d-$closed  $2q-$form $\Omega_{s}$ inducing $\omega -\partial\alpha - \bar\partial\beta$ on $X_{s}$. Then, the $(q, q)-$part of $\Omega_{s}$  is $\omega_{\vert X_{s}}$ so it is  strongly Lelong positive on $X_{s}$ and then also on  $\pi^{-1}(\Delta_{s})$, if $\Delta_{s}$ is small enough around $s$, because the strong Lelong positivity is an open property (and in our $\mathscr{C}^{\infty}-$trivialisation used to construct $\Omega_{s}$ the complex structure varies smoothly). The theorem above gives then the properness on $\Delta_{s}$ of the connected components of the space $\mathcal{C}_{q}(\pi_{\vert \Delta_{s}})$.\\
 To conclude, consider now a connected component $\Gamma$ of $\mathcal{C}_{q}(\pi)$ and remark that two $q-$cycles in $\Gamma$ which are contained in $\pi^{-1}(\Delta_{s})$ have not only the same image in $H_{2q}(\mathcal{X}, \mathbb{C})$ but also in $H_{2q}(\pi^{-1}(\Delta_{s}), \mathbb{C})$ because the natural map
  $$H_{2q}(\pi^{-1}(\Delta_{s}), \mathbb{C}) \to H_{2q}(\mathcal{X}, \mathbb{C}) $$
  is injective (in fact bijective). This implies that the integral of $\Omega_{s}$ on any cycle in   $\Gamma \cap \mathcal{C}_{q}(\pi^{-1}(\Delta_{s}))$ is constant, so $\Gamma \cap \mathcal{C}_{q}(\pi^{-1}(\Delta_{s}))$  has only finitely many connected components, each of them being proper on $\Delta_{s}$.$\hfill \blacksquare$
  
  \parag{Terminology} In the case of a compact complex connected manifold $X$ of dimension $n$ a {\bf strongly Gauduchon} form in the sense of  [B.15] is a smooth $d-$closed $(2n-2)-$form $\omega$ on $X$ such that its $(n-1, n-1)-$part is strongly positive in the sense of Lelong (in this case it is the ``usual'' sense). Assuming that $X$ satisfies the $\partial\bar\partial-$lemma, any Gauduchon form\footnote{so a $(n-1, n-1)$ smooth differential form  positive definite and $\partial\bar\partial-$closed.} gives rise to a strongly Gauduchon form via the method described above. And, thanks to [G.77] a Gauduchon form always exists on a compact complex  connected manifold.

 \subsection{Estimation of volumes}
 
  Let now consider a holomorphic family $\pi : \mathcal{X} \to  D$ of compact complex connected manifolds of dimension $n$ parametrized by the unit disc in $\C$. So $\mathcal{X}$ is a smooth complex manifold of dimension $n + 1$, and  we fix a smooth relative Gauduchon form $\omega$ on $\mathcal{X}$. Then $\omega$ is a smooth $\pi-$relative differential form on $\mathcal{X}$ on type $(n-1, n-1)$ which is positive definite in the fibres and satisfies $\partial_{/\pi}\bar\partial_{/\pi} \omega \equiv 0$ on $\mathcal{X}$. \\
   Using Ehresmann's theorem  we may assume that we have a $\mathscr{C}^{\infty}$ trivialization
  $$\xymatrix{ D \times X_{0} \ar[rr]^{\Phi} \ar[rd]_{pr} &\quad & \mathcal{X} \ar[ld]^{\pi}\\ \quad & D & \quad} $$
  of the map $\pi$ where $X_{0} := \pi^{-1}(0)$. Then the form $\Phi^{*}(\omega)$ is a smooth differential form on $X_{0}$ of degree $2n-2$  depending in a smooth way of the real variables $x$ and $y$ in $D$, and we then identify $\omega$ with the smooth (absolute) $(2n-2)-$differential form of $(\Phi^{-1})^{*}(\Phi^{*}(\omega))$ on $\mathcal{X}$.
  
   We shall define
  $$ \frac{\partial \omega}{\partial x} := (\Phi^{-1})^{*}(\frac{\partial (\Phi^{*}(\omega)}{\partial x}) \quad {\rm and } \quad \frac{\partial \omega}{\partial y} := (\Phi^{-1})^{*}(\frac{\partial (\Phi^{*}(\omega)}{\partial y}) $$
  on $\mathcal{X}$. Then $\frac{\partial \omega}{\partial x}$ and $\frac{\partial \omega}{\partial y}$ are smooth $(2n-2)-$differential forms on $\mathcal{X}$ and satisfy
  \begin{equation*}
   d\omega = \frac{\partial \omega}{\partial x}.\pi^{*}(dx) + \frac{\partial \omega}{\partial y}.\pi^{*}(dy) + d_{/\pi}\omega  . \tag{0}
   \end{equation*}
   where $d_{/\pi}$ is the $\pi-$relative differential.\\
  Now define the $\pi-$relative smooth differential forms  $\nabla_{x}\,\omega$ and $\nabla_{y}\,\omega$ on $\mathcal{X}$  as the $(n-1, n-1)-$parts of the restrictions of $\frac{\partial \omega}{\partial x}$ and $\frac{\partial \omega}{\partial y}$ to each fibre of $\pi$. We may also consider $\nabla_{x}\,\omega$ and $\nabla_{y}\,\omega$ as smooth hermitian forms on the holomorphic vector bundle $\Lambda^{n-1}(T_{\mathcal{X}/D})$ on $\mathcal{X}$. Then the following lemma is given by an elementary compactness argument, because $\omega$ is a smooth positive definite hermitian form on this vector bundle.
  
  \begin{lemma}\label{simple}
  For any compact set $K$ in $D$ there exists a constant $C_{K} > 0$ such on $\pi^{-1}(K)$ we have the following inequalities 
  \begin{equation*}
   \vert \nabla_{x}\,\omega\vert \leq C_{K}.\omega \quad {\rm  and} \quad \vert \nabla_{y}\,\omega\vert \leq C_{K}.\omega  \tag{1}
   \end{equation*}
   between hermitain forms on the vector bundle $\Lambda^{n-1}(T_{\mathcal{X}/D})$ on $\pi^{-1}(K)$.$\hfill \blacksquare$
   \end{lemma}
   
   The next lemma will allow us to prove our main estimate of the integral of $\omega$ over an analytic family of relative divisors in $\mathcal{X}$.
   
\begin{lemma}\label{crucial}
  Consider over an open disc $D' \subset D$ a closed and reduced complex hypersurface $\mathcal{H}$ in $\pi^{-1}(D')$ which is $(n-1)-$equidimensional on $D'$. Then, for any smooth differential $1-$form 
  $\varphi := \varphi_{x}.dx + \varphi_{y}.dy$ on $D'$ we have at the generic point of $\mathcal{H}$ the equality
   \begin{equation*}
   d\omega\wedge \pi^{*}(\varphi)_{\vert \mathcal{H}} = \big(\nabla_{x}\,\omega \wedge\pi^{*}(\varphi_{y}.dx\wedge dy) - \nabla_{y}\,\omega \wedge \pi^{*}(\varphi_{x}.dx\wedge dy)\big)_{\vert \mathcal{H}} \tag{2}
   \end{equation*}
   \end{lemma}
   
   \parag{proof} Near a generic point in $\mathcal{H}$ we may assume that we have a local holomorphic trivialisation
   $$ \xymatrix{D''\times U \subset D''\times V \ar[rr]^{\simeq} \ar[rd]_{pr} & \quad & \mathcal{H} \subset \mathcal{X} \ar[ld]_{\pi}\\ \quad & D'' & \quad } $$
   where $D''$ is a small open disc in $D$ and $V = U\times \Delta$  is the product of an open polydisc in $\C^{n-1}$ by an open disc $\Delta$ in $\C$. On a fibre of the restriction of $\pi$ to $\mathcal{H}$  the form $d\omega\wedge \pi^{*}(\varphi)$ is of pure type $(n, n)$. \\
   Remark that $d_{/\pi}\omega \wedge \pi^{*}(dx)$ and $d_{/\pi}\omega \wedge \pi^{*}(dy)$ vanish on $\mathcal{H} \simeq D''\times U$ because such  forms present only the types $(n, n-1)$ and $ (n-1, n)$ in the variables $t_{1}, \dots, t_{n-1}$ in $U \subset \C^{n-1}$. Then, using $(0)$ and the fact that  the restrictions to $\mathcal{H}$ of the forms 
   $$  \frac{\partial \omega}{\partial x}.\pi^{*}(\varphi_{y}.dx\wedge dy) \quad {\rm and} \quad  \frac{\partial \omega}{\partial y}.\pi^{*}(\varphi_{x}.dx\wedge dy) $$
   come only from the $(n-1, n-1)-$parts of the restriction to the fibre of $\pi$ of $ \frac{\partial \omega}{\partial x}$ and $ \frac{\partial \omega}{\partial y}$ respectively because the restriction to $\mathcal{H}$ of $\pi^{*}(dx \wedge dy)$ is of type $(1, 1)$, the conclusion follows, by definition of $\nabla_{x}\,\omega$ and $\nabla_{y}\,\omega$.$\hfill \blacksquare$\\
   
   The previous lemma gives the following simple estimate on $\mathcal{H}$.
   
   \begin{cor}\label{obvious}
   We keep the notations of the previous lemma. Assuming that the $1-$form $\varphi$ has a compact support in $D'$ contained in the compact set $K$ in $D$ we have on $\mathcal{H}$ the inequality
   \begin{equation*}
    \vert \big(d\varphi \wedge \pi^{*}(\varphi)\big)_{\mathcal{H}} \vert \leq C_{K}.(\vert \varphi_{x}\vert + \vert \varphi_{y}\vert).(\omega\wedge\pi^{*}(dx\wedge dy))_{\vert \mathcal{H}}. \tag{3}
    \end{equation*}
    between (smooth) real $2n-$differential forms on $\mathcal{H}$ where $C_{K}$ is the constant introduced in the lemma \ref{simple}.$\hfill \blacksquare$\\
   \end{cor}
   
   The following classical  lemma will be useful to conlude our estimation of the volume.
   
    \begin{lemma}\label{conj.}
 Let $U$ a bounded domain in $\C$. Assume that there is a point $s_{0} \in U$ such that any point in $U$ can be joined  to $s_{0}$ by a $\mathscr{C}^{1}$ path in $U$  with length uniformly bounded by a number $L$. For instance, we can restrict ourself to $U := D^{*}\setminus A$ where $D^{*}$ is a punctured open disc in $\C$ and $A$ a closed discrete subset in $D^{*}$. Let $\eta : U \to \mathbb{R}^{+*}$ be a $\mathscr{C}^{1}$  function on $U$ such that 
  \begin{equation*}
  \vert \frac{\partial \eta}{\partial x} \vert \leq  C.\eta \quad {\rm and } \quad \vert \frac{ \partial \eta}{\partial y} \vert \leq  C.\eta \quad {\rm    in} \ U \tag{@}
  \end{equation*} 
  where $C > 0$ is a given constant. Then this implies that the function $\eta$ is uniformly bounded on $U$ by $\eta(s_{0}).exp(2C.L)$.
   \end{lemma}
 
 \parag{Proof} Firstly consider the case of an open interval $]a, b[$ in $\mathbb{R}$ and let $L := b - a$. If $\eta : ]a, b[ \to \mathbb{R}^{+*}$ is a $\mathscr{C}^{1}$ function such that $\vert \eta'(x) \vert \leq C.\eta(x) \quad \forall x \in ]a, b[$; we have 
 $$ -C \leq \frac{\eta'(x)}{\eta(x)} \leq C $$
 which gives, after integration on $[x_{0}, x] \subset ]a, b[$
 $$ - C.\vert x - x_{0}\vert \leq Log \frac{\eta(x)}{\eta(x_{0})} \leq C.\vert x-x_{0}\vert $$
 and then, for any choice of $x_{0}\in ]a, b[$ the estimate $\eta(x) \leq \eta(x_{0}).\exp(C.L)$ is valid.\\
 In the case of the bounded domain $U \subset \C$, choose for any $s \in U$ a $\mathscr{C}^{1}$ path of length strictly less than $L$ and consider a parametrization of this path extended a little around $s_{0}$ and $s$ :
 \begin{align*} 
 & \Phi : ]a, b[ \to U \quad   x := \varphi(t), y := \psi(t) \quad {\rm with} \quad  \varphi'(t)^{2}+ \psi'(t)^{2} \equiv 1 \\
 & {\rm and} \quad  \Phi(a+\varepsilon) = s_{0} \quad {\rm and} \quad \Phi(b - \varepsilon) = s \quad {\rm with} \quad 0 < \varepsilon \ll 1 .
 \end{align*}
 Then apply the one dimensional case to the function  $F(t) := \eta(\Phi(t))$ which satisfies $\vert F'(t)\vert \leq 2C.F(t)$. This gives $\eta(s) \leq \eta(s_{0}).exp(2C.L)\quad \forall s \in U$.$\hfill \blacksquare$\\

   \subsection{proof of the theorem \ref{main}}
   
   \parag{First step} 
For each $s \in  D^{*}$ we know that $X_{s}$ admits a strongly Gauduchon form and then there exists a strongly Gauduchon form for the family over a small disc around $s$ (see section 3.1). So we have a smooth $d-$closed $2(n-1)-$form such that its $(n-1, n-1)$ part is positive definite on each fibre. This implies that each connected component of the relative cycle space $\mathcal{C}_{n-1}(\pi^{*})$ is proper on $ D^{*}$, where $\pi^{*}$ is the restriction of $\pi$ to $\pi^{-1}(D^{*})$:\\
in a connected component $\Gamma$ of the relative cycle space $\mathcal{C}_{n-1}(\pi^{*})$, for any disc $\Delta \subset D^{*} $ all cycles in  $\pi^{-1}(\Delta)\cap \Gamma$ are homologuous in $\pi^{-1}(\Delta)$ because we have a $\mathscr{C}^{\infty}-$ relative isomorphism $\mathcal{X} \simeq D \times X_{0}$. So if we dispose of a $d-$closed $2(n-1)-$ smooth form $\Omega$ on $\pi^{-1}(\Delta)$ with a positive definite $(n-1, n-1)$ part on fibres we see that the volume for $(n-1)-$cycles in $\Gamma \cap \pi^{-1}(\Delta)$ relative to the ``volume'' defined by $\Omega$ has to be constant, proving the properness on $\Delta$ of $\Gamma \cap \pi^{-1}(\Delta)$.

\parag{Second step} Applying the  theorem 1.0.2 of [B.15] we obtain that the minimum of the algebraic dimension of the fibres over $D^{*}$ is obtained on a dense (in fact the complement of a countable subset) subset of $D^{*}$. This implies  that it is enough, under the hypothesis of the theorem \ref{main}, to prove that the conclusion holds on a small open disc $D_{0}$ with center $0$ in $D$ because this ``weak'' version of this theorem may be apply then to a small disc centered at  any point point in $D^{*}$ and this gives the theorem \ref{main}.\\

Thanks to the second step, we may now assume that there exists a positive constant $C$ such that the inequalities of the the lemma \ref{simple} holds on $\mathcal{X}$.\\
Now, using the proposition \ref{classique} (recalled from [B.15]) it is enough to show that each irreducible component $\Gamma$  of the space $\mathcal{C}_{n-1}(\pi)$ of relative $(n-1)-$cycles in $\mathcal{X}$ is proper over $D$.\\

\parag{Step 3} 
Let $\Gamma$ be an irreducible  component of $\mathcal{C}_{n-1}(\pi^{*})$ with reduced generic cycle and such that the projection  $p : \Gamma \to  D^{*}$ is surjective. As the map $p$ is proper, let $\tau : \pi^{-1}( D^{*}) \to S$ be a  Stein factorization (we may assumed that $S$ is normal) of $p$ and let $\tilde{\tau} : S \to  D^{*}$ the corresponding  proper finite surjective map.\\
Now let $\omega$ be a smooth $\pi-$relative Gauduchon form on $\mathcal{X}$ and define on $\Gamma$ the function
$$ \theta : \Gamma \to \mathbb{R}^{+*} \quad  \gamma \mapsto \int_{\gamma}\ \omega .$$
 As the function $\theta$ is continuous,  pluri-harmonic along the fibres of $p$ as a distribution\footnote{meaning that it is $\partial\bar\partial-$closed as a distribution along the fibres of  $p :\Gamma \to D^{*}$.} and proper, it is constant in the fibres of $\tau$ which are the connected components of fibres of $p$ over $D^{*}$. So it  defines a continuous function $\tilde{\theta} : S \to \mathbb{R}^{+*}$. This situation corresponds to the following diagram:

 $$\xymatrix{ \quad & \Gamma \ar[ldd]_{\theta} \ar[ddr]^{p}  \ar[dd]^{\tau}& G \ar[l]^{p_{1}} \ar[r]  \ar[rd]^{p_{2}}&  \Gamma\times \mathcal{X} \ar[d]^{p_{2}} \\
 \quad &\quad & \quad& \mathcal{X} \ar[d]^{\pi} \\
 \mathbb{R}^{+*} & S \ar[l]_{\tilde{\theta}} \ar[r]^{\tilde{\tau}} &\quad  D^{*} \ar[r] & D}$$
 
 \parag{Step 4} Let $\tilde{A}_{0} \subset S$ be the set of points $\sigma \in S$ where the fibre $\tau^{-1}(\sigma)$ is not contained in the union of the singular set of $\Gamma$ with the critical set of $p : \Gamma \to D^{*}$. Then $\tilde{A}_{0}$ is a closed analytic subset in $S$ with no interior point. So $A_{0} := \tilde{\tau}(A_{0})$ is a closed discrete subset in $D^{*}$. Let $A_{1} \subset D^{*}$ be the ramification set of $\tilde{\tau}$. It is also closed and discrete in $D^{*}$ and finally put $A := A_{0}\cup A_{1}$.\\
 We shall show now that the continuous function  
 $$ \eta : D^{*} \to \mathbb{R}^{+*}, s \mapsto \eta(s) := Trace_{\tilde{\tau}}(\tilde{\theta})(s) $$
 satisfies the following properties :
 \begin{enumerate}[i)]
 \item The partial derivatives $\frac{\partial \eta}{\partial x}$ and $\frac{\bar\partial \eta}{\partial y}$ are continuous on $D^{*}\setminus A$ (with  $s := x + i.y$).
 \item They satisfy the inequalities $\vert \frac{\partial \eta}{\partial x} \vert \leq C.\eta$ and $\vert \frac{\partial \eta}{\bar\partial y} \vert \leq C.\eta$ on $D^{*}\setminus A$, where $C$ is the constant introduced above (see lemma \ref{simple} and the end of step 2).
 \end{enumerate}
 Note that it is enough to prove the properties i) and ii) near each point in $D^{*}\setminus A$.\\
 Choose a point $s_{0}\in D^{*}\setminus A$ and fix an open disc $\Delta_{0} \subset D^{*}\setminus A$ with center $s_{0}$. Choose now in $p^{-1}(s_{0})$  generic point $\gamma_{0}$ and, as $\gamma_{0}$ is a smooth point of $\Gamma$, such that $p$ has rank $1$ at this point, we can find a smooth locally closed curve $\Sigma$ through $\gamma_{0}$ in $\Gamma$ such that the restriction of $p$ to $\Sigma$ induces an isomorphism of $\Sigma$ onto $\Delta_{1} $ where $\Delta_{1} \subset \Delta_{0}$ is an open disc with center $s_{0}$.

 Define, for $s \in \Delta_{1}$ the cycle $Y_{s}$ as the cycle corresponding to the point $p_{\vert \Sigma}^{-1}(s)$ in $\Sigma \subset \Gamma$. This defines an analytic family of $\pi-$relative $(n-1)-$cycles. So the function $\eta_{1}$ which is given by $s \mapsto \int_{Y_{s}} \ \omega$ is continuous on $\Delta_{1}$ (see the proposition IV  2.3.1  in [B-M 1]) and coincides with $\frac{1}{\delta}.\eta$ where $\delta$ is the degree of the map $\tilde{\tau}$. We shall compute now the partial derivatives in the distribution sense of this function $\eta_{1}$.\\

 Let $\mathcal{Y} \subset \Delta_{1}\times \mathcal{X}$ be the graph of the analytic family $(Y_{s})_{s \in \Delta_{1}}$ and let $\varphi$ be in $\mathscr{C}_{c}^{\infty}(\Delta_{1})$. Then we have
 \begin{align*}
&  \langle d\eta_{1}, \varphi.dy \rangle = -  \langle \eta_{1}, \frac{\partial \varphi}{\partial x}.dx \wedge dy \rangle = - \int_{\mathcal{Y}} \  \omega\wedge \pi^{*}(\frac{\partial \varphi}{\partial x}.dx\wedge dy) \\
& \qquad =  - \int_{\mathcal{Y}}  \  \omega\wedge \pi^{*}(d\varphi\wedge dy) = - \int_{\mathcal{Y}}  \  \omega\wedge d(\pi^{*}(\varphi. dy))\quad {\rm then} \\
&\langle d\eta_{1}, \varphi.dy \rangle =  \int_{\mathcal{Y}} \ \frac{\partial \omega}{\partial x}\wedge \pi^{*}(\varphi.dx\wedge dy) = \int_{\mathcal{Y}} \nabla_{x}\,\omega \wedge \pi^{*}(\varphi.dx\wedge dy) 
\end{align*}
by Stokes formula and using the lemma \ref{crucial}.

We conclude that the distribution $\frac{\partial \eta_{1}}{\partial x}$ is equal to the continuous function
$$s \mapsto \int_{Y_{s}} \ \nabla_{x}\,\omega$$
with the estimate 
$$ \vert \frac{\partial \eta_{1}}{\partial x} \vert \leq C.\eta_{1}$$
deduced from the inequality of the lemma \ref{simple}.\\

We have the analogous result for  $\frac{\partial \eta_{1}}{\partial y}$.\\

So the lemma \ref{conj.} applies and the function $\theta$ on $\Gamma$ is uniformly bounded. Then the closure of $\Gamma$ in $\mathcal{C}_{n-1}(\mathcal{\pi})$ is proper over $D$, using again the characterization of compact subsets in $\mathcal{C}_{n-1}(\mathcal{X})$ (see [B-M 1] IV  th.2.7.20) .\\

\parag{Final step} Consider now any irreducible component $\tilde{\Gamma}$ of $\mathcal{C}_{n-1}(\pi)$. Remark that to prove the properness of $\tilde{\Gamma}$ on $D$ it is enough to consider the case where the generic cycle of $\tilde{\Gamma}$ is reduced.\\
Then either $\tilde{\Gamma}$ is contained in some  $\mathcal{C}_{n-1}(X_{s})$ for some $s \in D$, and then it is a connected closed subset in $\mathcal{C}_{n-1}(X_{s})$ and then it is compact, or  $\tilde{\Gamma} \setminus p^{-1}(0)$ is a closed irreducible subset in $\mathcal{C}_{n-1}(\pi^{*})$ which is proper and surjective on $D^{*}$ because we may apply the previous result. In this cases we conclude also that $\tilde{\Gamma}$ (which is the closure of $\tilde{\Gamma} \setminus p^{-1}(0)$) is proper over $D$.\\
As any irreducible component of $\mathcal{C}_{n-1}(\pi)$ is proper over $D$, the conclusion follows.$\hfill \blacksquare$\\

 For the convenience of the reader let me recall the proposition 4.0.8 of [B.15] which is an easy generalization of an old result of F. Campana (see [C.81]).
 
 \begin{prop}\label{classique}
 Let $\pi : \mathcal{X} \to S$ be a proper surjective holomorphic $n-$equidimensional map between two irreducible complex spaces. Assume that any irreducible component of the complex space $\mathcal{C}_{n-1}(\pi)$ is proper over $S$. Then there exists a countable union $\Sigma$ of closed irreducible analytic subsets with no interior points in $S$ and a non negative integer $a$ such that:
 \begin{enumerate}
 \item For any $s \in S \setminus \Sigma$ the algebraic dimension of $X_{s}$ is equal to $a$;
 \item For all $s \in S$ the algebraic dimension of $X_{s}$ is at least equal to $a$.
 \end{enumerate}
 \end{prop}

\bigskip
   
    We conclude by noticing that there exists an analytic family of smooth complex compact surfaces of the class VII (not K\"{a}hler) parametrized by a disc Æ such that the central fibre has algebraic dimension 0 and all other fibres have algebraic dimension equal to $1$, see [F-P.09].\\
This shows that in our theorem \ref{main} some ``K\"{a}hler type'' assumption on the general fibre $X_{s}$ cannot be avoided in order that the ÒgeneralÓ algebraic dimension gives a lower bound for the algebraic dimensions of all fibres. \\
Note that our assumption that the generic fibres satisfy the $\partial\bar \partial-$lemma (in fact for the type $(n, n-1)$) is a rather weak such assumption.

\newpage

 \section{Bibliography}
 
 \bigskip
 
 \bigskip
 
 \bigskip
 
 \bigskip
 
 \begin{enumerate}
 \item{[B.15]}  Barlet, D. {\it Two semi-continuity results for the algebraic dimension of compact complex manifolds},  J. Math. Sci. Univ. Tokyo 22 (2015), no. 1, pp.39-54.
 \item{[B-M 1]} Barlet, D.  et  Magn\'usson, J.  {\it Cycles analytiques complexes. I. Th\'eor\`{e}mes de pr\'eparation des cycles}, Cours Sp\'ecialis\'es  22., Soci\'et\'e Math\'ematique de France, Paris, 2014. 525 pp.
 \item{[Bi.64]} Bishop, E. {\it Conditions for the analyticity of certain sets}, Mich. Math. Journal (1964) p.289-304.
 \item{[C.82]} Campana, F., {\it Sur les diviseurs non-polaires d'un espace analytique \\ 
 compact}, J. Reine Angew. Math. 332 (1982), 126-133.
 \item{[Fu.82]} Fujiki, A., {\it Projectivity of the space of divisors of a normal compact complex space}, Publ. Res. Inst. Math. Sci. 18 (1982), no. 3, 1163-1173.
  \item{[F-P.09]} Fujiki, A. and Pontecorvo, M. {\it Non-upper continuity of algebraic dimension for families of compact complex manifolds} \\
 arXiv: 0903.4232v2 [math. AG] .
 \item{[G.77]} Gauduchon, P., {\it Le th\'eor\`eme de l'excentricit\'e nulle}, C. R. Acad. Sci. Paris, S«er. A 285 (1977), 387-390.
 \end{enumerate}

 \end{document}